\documentclass[12pt]{amsart}
\usepackage{verbatim}
\newcommand{\car}{\operatorname{char}}
\newcommand{\eff}{{\operatorname{eff}}}
\newcommand{\rat}{{\operatorname{rat}}}
\newcommand{\num}{{\operatorname{num}}}
\renewcommand{\o}{{\operatorname{o}}}
\newcommand{\Mot}{\operatorname{Mot}}
\newcommand{\Ker}{\operatorname{Ker}}

\newcommand{\inj}{\hookrightarrow}
\newcommand{\iso}{\stackrel{\sim}{\longrightarrow}}

\newcommand{\un}{{\mathbf{1}}}

\newcommand{\F}{\mathbf{F}}
\renewcommand{\P}{\mathbf{P}}
\newcommand{\Q}{\mathbf{Q}}
\newcommand{\Z}{\mathbf{Z}}

\newcommand{\cf}{{\it cf. }}

\newcommand{\prf}{\noindent{\bf Proof. }}

\newtheorem{thm}{Theorem}
\newtheorem{lemma}{Lemma}
\newtheorem{prop}{Proposition}
\newtheorem{cor}{Corollary}
\newtheorem{conj}{Conjecture}
\theoremstyle{definition}
\newtheorem{defn}{Definition}
\newtheorem{rk}{Remark}
\newtheorem{rks}[rk]{Remarks}
\newtheorem{ex}{Example}

\begin{document}

\title{Number of points of function fields over finite fields}
\author{Bruno Kahn}
\address{Institut de Math\'ematiques de Jussieu\\175--179 rue du 
Chevaleret\\75013 Paris\\France}
\email{kahn@math.jussieu.fr}
\date{October 28, 2002}
\maketitle

\enlargethispage*{20pt}

\section*{Introduction}
Let $k$ be a field; if $\sim$ is an adequate equivalence relation on
algebraic  cycles, we denote by $\Mot_{\sim}(k)$ or simply $\Mot_{\sim}$
the category of  motives modulo $\sim$ with rational coefficients, and by
$\Mot_{\sim}^\eff$ its  full subcategory consisting of effective motives
\cite{scholl}\footnote{With  notation as in \protect\cite{scholl}, an
object of $\Mot_\sim$ is 
\emph{effective} if it is isomorphic to an triple $(X,p,n)$ with $n\le
0$; one  may then find such a triple with $n=0$.}. We use the convention
that the functor $X\mapsto h(X)$ from smooth projective $k$-varieties to
$\Mot^\eff_\sim$ is covariant. We shall in  fact only consider the two
extreme cases: rational equivalence (rat) and  numerical equivalence
(num).

Using the point of view of birational motives (developed
jointly with Sujatha \cite{ks}), we give a proof almost without
cohomology (see proof of Lemma \ref{l1}) of a recent result of Esnault on
the existence of rational points for smooth projective varieties with
``trivial" Chow group of zero-cycles over a finite field \cite{ln}. We
also prove that the number of rational points modulo $q$ is a stable
birational invariant of smooth projective varieties over
$\F_q$:  the
idea of considering effective motives and their divisibility by the
Lefschetz motive was anticipated by Serre \cite{serre}. This answers a question of Koll\'ar; the
$3$-dimensional case had been dealt with by Lachaud and Perret earlier
\cite{lp}. However, as was pointed out by Chambert-Loir, this birational
invariance in fact follows from much earlier work of Ekedahl
\cite{ek}, who does not use any form of resolution of singularities!

\section{Birational motives}

\begin{defn} The category $\Mot_{\sim}^\o$ is the Karoubian envelope (or 
idempotent completion) of the quotient of $\Mot_\sim^\eff$ by the ideal
$J$  consisting of morphisms factoring through an object of the form
$M\otimes L$, where $L$ is the Lefschetz motive. This is a tensor
additive category. If $M\in 
\Mot_\sim^\eff$, we denote by $\bar M$ its image in $\Mot^\o_\sim$.
\end{defn} 

\begin{lemma}[\protect{\cite[Lemmas 5.3 and 5.4]{ks}}] \label{l0} Let
$X,Y$ be  two smooth 
projective irreducible $k$-varieties. Then, in $\Mot_\rat^\o$, we have
\[Hom(\bar h(X),\bar h(Y))=CH_0(Y_{k(X)})\otimes\Q.\]
\end{lemma}

(Let us briefly recall the proof: for $X,Y$ smooth projective, let
$I(X,Y)$ be the subgroup of $CH^{\dim Y}(X\times Y)\otimes\Q$ formed of
those correspondences which vanish on $U\times Y$ for some dense open
subset $U$ of $X$. Then $I$ is an ideal in the category of rational Chow
correspondences: the proof \cite[Lemma 5.3]{ks} is a slight
generalisation of the argument in \cite[Ex. 16.1.11]{fulton}. It is even
monoidal, and extends to a monoidal ideal $I$ in $\Mot_\rat^\eff$, which
obviously contains $J$. Using de Jong's theorem \cite[Th. 4.1]{dJ} and
Chow's moving lemma, one sees that $I\otimes\Q=J\otimes\Q$ \cite[Lemma
5.4]{ks}. In characteristic $0$, one may remove the coefficients $\Q$ by
using Hironaka's resolution of singularities.)

\begin{ex}\label{e1} Let $X$ be smooth and projective over $k$. Then $\bar
h(X)\simeq \un$ in $\Mot_\rat^\o$ if and only if
$CH_0(X_{k(X)})\otimes\Q\simeq \Q$ (write $h(X)\simeq \un\oplus
h(X)_{\ge 1}$ in $\Mot_\rat^\eff$).
\end{ex}

\begin{rk}\label{r1} If $K$ is the function
field of a smooth  projective variety $X$, we may define a motive $\bar
h(K)\in \Mot_\rat^\o$ as follows.  If $Y$ is another smooth projective
model of $K$, then [the closure of] the graph of a birational isomorphism
from $X$ to $Y$ defines an isomorphism $\bar h(X)\iso \bar h(Y)$. If
there is a third model $Z$, then the system of these isomorphisms is
transitive, so defining $\bar h(K)$ as the direct limit of the $\bar
h(X)$ for this type of isomorphisms makes sense and is (canonically)
isomorphic to any of the
$\bar h(X)$. This construction is functorial for inclusions of
fields. If $\car k=0$, it is even functorial for
$k$-places by \cite[Lemma 5.6]{ks}, although we
won't use this. (Extending it to arbitrary function fields in
characteristic $p$ would demand more work.)

Note that if $K\subseteq L$, then $\bar h(K)$ is a direct summand of
$\bar h(L)$: to see this, consider smooth projective models $X,Y$ of $K$
and $L$. Let $f:U\to X$ be a corresponding dominant morphism, where $U$ is
an open subset of $Y$. By Noether's normalisation theorem, we may find an
affine open subset $V\subset U$ such that the restriction of $f$ to $V$
factors through $X\times \P^n$, where $n=\dim Y-\dim X$. Since $\bar
h(X\times\P^n)\iso \bar h(X)$, we are reduced to the case where $L/K$ is
finite, and then it follows from a transfer argument. In particular,
$\bar h(X)\simeq\un$ if
$X$ is unirational, as expected in
\cite{serre}. The converse is not true: an Enriques surface $X$ verifies
$\bar h(X)\simeq \un$ by \cite{bkl} (see example \ref{e1}), but is not
rational, hence not unirational over a field of characteristic $0$
because
$Pic(X)$ contains a $\Z/2$ summand (this counterexample was explained by
Jean-Louis Colliot-Th\'el\`ene).
\end{rk}

For $\sim=\num$, the category $\Mot_{\sim}^\o$ is abelian semi-simple 
\cite{jannsen}. From \cite[Prop. 2.1.7]{ak}, we therefore get:

\begin{prop}\label{p0} a) The projection functor $\pi:\Mot_\num^\eff\to 
\Mot_\num^\o$ is essentially surjective (i.e. taking the karoubian envelope is 
irrelevant in the definition of $\Mot_\num^\o$).\\
b) $\pi$ has a section $i$ which is also a left and right adjoint.\\
c) The category $\Mot_\num^\eff$ is the coproduct of $\Mot_\num^\eff\otimes L$ 
and $i(\Mot_\num^\o)$, i.e. any object of  $\Mot_\num^\eff$ can be uniquely 
written as a direct sum of objects of these two subcategories.\\
d) The sequence
\[0\to K_0(\Mot_\num^\eff)\stackrel{\cdot L}{\longrightarrow} 
K_0(\Mot_\num^\eff)\to 
K_0(\Mot_\num^\o)\to 0\]
is split exact.\qed
\end{prop}

(In d), the injectivity on the left corresponds to the fact that the functor 
$-\otimes L$ is fully faithful.)

\begin{rk}\label{r2} a) In $\Mot_\num^\o$, we can extend the end of
Remark \ref{r1} as follows: let $K,L$ two function fields of smooth
projective varieties such that
$K\inj L(t_1,\dots,t_m)$ and $L\inj K(t_1,\dots, t_n)$ for some $m,n$.
Then $\bar h(K)\simeq \bar h(L)$. To get such a result in $\Mot_\rat^\o$,
one would need to have enough information on the algebra $End(\bar h(K))$.

b) Proposition \ref{p0} b) shows via Remark \ref{r1}
that to  any function field $K/k$ one may canonically associate an
effective numerical  motive $h(K)\allowbreak\in \Mot_\num^\eff$, which is
a direct summand of $h(X)$ for any  smooth projective model $X$ of $K$ (if
any).
\end{rk}

The following conjecture was suggested by Luca Barbieri-Viale:

\begin{conj}[\cf \protect{\cite[Conj. 0.0.11]{voebei}}] For any field $k$, the 
projection functor $\Mot_\rat^\eff\to \Mot_\rat^\o$ has a right adjoint.
\end{conj}

\enlargethispage*{20pt}

\section{Number of rational points modulo $q$}

From now on, $k=\F_q$ is a field with $q$ 
elements. Then, for all $n\ge 1$, the assignment
\[X\mapsto |X(\F_{q^n})|=\deg(\Delta_X\cdot F_X^n)\]
for a smooth projective variety $X$, where $\Delta_X$ is the class of the 
diagonal and $F_X$ is the Frobenius endomorphism viewed as a correspondence, 
extends to a ring homomorphism
\begin{equation}\label{eq1}
\#_n:K_0(\Mot^\eff_\num)\to\Q
\end{equation}
by the rule $\#_n(X,p)= \deg({}^tp\cdot F_X^n)$ if $p=p^2\in End (h(X))$, \cf 
\cite[p. 80]{kleiman}.

\begin{lemma}\label{l1} The homomorphisms $\#_n$ take their values in $\Z$.
\end{lemma}

\prf It is enough to prove this for $n=1$. More conceptually, we have 
$\deg({}^tp\cdot F_X)=Tr(p\circ F_X)$ in the rigid tensor category
$\Mot_\rat$.  We may compute this trace after applying a Weil cohomology
$H$, e.g. $l$-adic  cohomology. (We then have to consider $H(X)$ as a
$\Z/2$-graded vector space  and compute a super-trace.) Let $H(X)=V\oplus
W$, with $V=\Ker(H(p)-1)$ and 
$W=\Ker(H(p))$. Since $F_X$ is a central correspondence, it commutes with 
$p$, hence $H(F_X)$ respects $V$ and $W$ and
\begin{multline*}
Tr(p F_X) = Tr(H(pF_X)) =
Tr(H((pF_X)_{|V})+Tr(H((pF_X)_{|W})\\
=Tr(H(F_X)_{|V}).
\end{multline*}

\begin{sloppypar}
Since the minimum polynomial of $H(F_X)$ kills $H((F_X)_{|V})$, the
eigenvalues  of the latter are algebraic integers. Hence $Tr(p F_X)$ is
an algebraic integer  and therefore is in $\Z$.\qed
\end{sloppypar}

\enlargethispage*{20pt}

\begin{thm}\label{t1} The homomorphism \eqref{eq1} induces a ring
homomorphism
\[\overline\#_n:K_0(\Mot^\o_\num)\to\Z/q^n.\]
\end{thm}

This follows from Lemma \ref{l1} and Proposition \ref{p0} d) (note that 
$\#_n(L)=q^n$).\qed

\begin{cor}[Esnault \protect{\cite{ln}}]\label{c1} Let $X$ be a smooth
projective variety  over $\F_q$ such that $CH_0(X_{\F_q(X)})\otimes\Q
=\Q$. Then
$|X(\F_q)|\equiv  1\pmod{q}$.
\end{cor}

\prf By Example \ref{e1}, one has $\bar h(X)\simeq \un$
in $\Mot_\rat^\o$, hence a fortiori in $\Mot_\num^\o$.\qed

\begin{cor}[\cf \protect{\cite[Th. 4]{ek}, \cite{lp}}]\label{c2} The
number of rational points modulo $q$ is a stable birational invariant of
smooth projective
$\F_q$-varieties.
\end{cor}

Indeed, two stably birationally isomorphic varieties have isomorphic motives in 
$\Mot_\rat^\o$. \qed

\begin{rks} a) Using Remark \ref{r2} a) we could strengthen Corollary
\ref{c2} as follows: for two smooth projective varieties $X,Y$,
$|X(\F_q)|\equiv |Y(\F_q)|\pmod{q}$ provided there exist
$m,n$ and dominant rational maps $X\times \P^m\to Y$, $Y\times\P^n\to
X$. However, this also follows from \cite{ek}.\\   
b) Using Remark \ref{r2} b) we may canonically  associate to any
function field $K/\F_q$  a series of integers $(a_n)_{n\ge 1}$  such
that, for all
$n$, $\overline \#_n(\bar h(K))\allowbreak = a_n\pmod{q^n}$  (see Remark
\ref{r1} for the definition of $\bar h(K)$). Naturally $a_n$ is not 
necessarily positive in general. More conceptually, we may associate to
$K$ its 
\emph{zeta function}, defined as the zeta function of the motive $i(\bar 
h(K))$.\\ 
c) Killing $L^\kappa$ instead of $L$ would yield congruences modulo 
$q^\kappa$ rather than modulo $q$, \cf \cite[\S3]{ln}; compare also 
\cite{voesli}. But one would lose the fact that function fields have motives as 
in the previous remarks.\\
d) Unfortunately the proof of Lemma \ref{l1} uses cohomology, hence the
proof of Theorem \ref{t1} is not completely cohomology-free.
\end{rks}

\section{A conjectural converse}\label{tate}

Note that the functions $\#_n$ of \eqref{eq1} extend to ring homomorphisms 
$K_0(\Mot_\num)\to \Z[1/q]$, still denoted by $\#_n$.

\begin{sloppypar}
\begin{thm}\label{t2} Assume that the Tate conjecture holds. Let $M\in
K_0(\Mot_\num)$ be  such that $\#_n(M)\in\Z$ for all $n\ge 1$. Then $M\in
K_0(\Mot_\num^\eff)$.  Conversely, if this implication holds for any
$M\in K_0(\Mot_\num)$, then the  Tate conjecture holds.
\end{thm}
\end{sloppypar}

\prf Write $M=\sum m_i[S_i]$, where $m_i\in\Z\setminus\{0\}$ and the $S_i$ are 
simple pairwise non-isomorphic motives. For each $i$, let $w_i$ be a Weil 
number of $S_i$, that is, a root of the minimum polynomial of $F_{S_i}$, and 
$K_i=\Q(w_i)$. Then $\#_n(M)=\sum m_i Tr_{K_i/\Q}(w_i^n)$. It follows from the 
assumption and from \cite[Lemma 2.8]{kleiman2} that $w_i$ is an algebraic 
integer for all $i$. (To apply {\it loc. cit.}, compute in a Galois extension of 
$\Q$ containing all $K_i$ and observe that the Tate conjecture implies that no 
$w_i$ is equal to a conjugate of $w_j$ for $i\ne j$ \cite[proof of Prop. 
2.6]{milne}.)

By Honda's theorem, for each $i$ there is an abelian variety $A_i$ over $\F_q$ 
and a simple direct summand $T_i$ of $h(A_i)$ whose Weil
numbers are the Galois  orbit of $w_i$ ({\it ibid}.). Reapplying Tate's
conjecture, we get that 
$S_i\simeq T_i$, hence $S_i$ is effective for all $i$.

To prove the converse, let $M\in \Mot_\num$ be simple and such that $F_M=1$.
Then $\#_n(M)=1$ for all $n$. Therefore $M$ is effective. Writing $M$ as
$(X,p)$ for  some smooth projective variety $X$, we have that $M$ is a direct
summand of $h^0(X)$ for weight reasons. It follows easily that $M\simeq \un$.
By \cite[Th.  2.7]{geisser}, this implies the Tate conjecture.\qed

\enlargethispage*{20pt}

\begin{sloppypar}
\begin{rk} Unfortunately we have to apply the Tate conjecture to
$S_i\otimes  T_i^*$ in the proof, hence cannot provide a hypothesis only
\text{involving 
$M$.}
\end{rk}
\end{sloppypar}

\begin{cor} The ring homomorphism
\[(\overline\#_n)_{n\ge 1}: K_0(\Mot_\num^\o)\to \prod_{n=1}^\infty
\Z/q^n\] 
is injective if and only if the Tate conjecture is true.\qed
\end{cor}

\subsection*{Acknowledgements} Part of this work was done during a
stay at the Newton Institute of  Mathematical Sciences in
September-October 2002; I thank it for its hospitality. I am grateful to
Dan Abramovich, Yves Andr\'e, Antoine Chambert-Loir, Jean-Louis
Colliot-Th\'el\`ene and Gilles Lachaud for helpful remarks. Especially
Colliot-Th\'el\`ene and Chambert-Loir pointed out the work of Lachaud and
Perret \cite{lp}, Cham\-bert-Loir pointed out the work of Ekedahl
\cite{ek} and Lachaud kindly communicated the letters of Serre
\cite{serre}. Andr\'e pointed out that Lemma \ref{l1} needed a proof,
helped with it and let me discover a mistake in the very first version of
this paper. Finally, the idea of this article was of course sparked by
H\'el\`ene Esnault's prior article
\cite{ln}.

\end{document}